\documentclass[10pt]{amsart}    
\title[Coeffective cohomology of symplectic aspherical manifolds ]
{Coeffective cohomology of symplectic aspherical manifolds}

\author{Hisashi Kasuya}
\usepackage{amssymb} 
\usepackage{amsmath} 
\usepackage{amscd}
\usepackage{amstext} 
\usepackage{amsfonts}
\usepackage[all]{xy}

\address[H.kasuya]{Graduate school of mathematical science university of tokyo japan }
\curraddr{}
\email{khsc@ms.u-tokyo.ac.jp}

\thanks{}

\newcommand{\C}{\mathbb{C}}
\newcommand{\R}{\mathbb{R}} 

\newcommand{\Ra}{\mathbb{Q}}
\newcommand{\Z}{\mathbb{Z}}
\newcommand{\g}{\frak{g}}

\theoremstyle{plain}
\newtheorem{theorem}{Theorem}[section] 
\theoremstyle{remark}
\newtheorem{remark}{Remark}
\theoremstyle{lemma}
\newtheorem{lemma}[theorem]{Lemma}
\theoremstyle{definition}
\newtheorem{definition}[theorem]{Definition}
\theoremstyle{proposition}
\newtheorem{proposition}[theorem]{Proposition}
\theoremstyle{corollary}
\newtheorem{corollary}[theorem]{Corollary}
\theoremstyle{remark}
\newtheorem{example}{Example}
\begin{document} 
\begin{abstract} 
We prove a generalization of the theorem which is  proved by Fernandez, Ibanez, and de Leon.
By this result, we give examples of non-K\"ahler manifolds which satisfy the property of compact K\"ahler manifolds concerning the coeffective cohomology.
\end{abstract}
\maketitle
\section{Introduction}
Let $(M,\omega)$ a compact $2n$-dimensional sympectic manifold.
Denote $A^{\ast}(M)$ the de Rham complex of $M$.
We call a differential  form $\alpha\in A^{\ast}(M)$ coeffective if $ \omega\wedge \alpha=0$, and  denote  the sub-DGA
\[A^{\ast}_{coE}(M)=\{\alpha\in A^{\ast}(M)\vert  \omega\wedge \alpha=0 \}.
\]
We call the cohomology $H^{\ast}(A^{\ast}_{coE}(M))$ the coeffective cohomology of $M$.
We also denote 
\[\tilde H^{\ast}(A^{\ast}(M))=\{[\alpha]\in H^{\ast}(A^{\ast}(M))\vert [\omega]\wedge [\alpha]=0\}.
\]
\begin{theorem}{\rm (\cite{Bou})}
Let $(M,\omega)$ be a compact K\"ahler manifold.
For $p\ge n+1$, we have an isomorphism
\[H^{p}( A^{\ast}_{coE}(M))\cong \tilde H^{p}(A^{\ast}(M)).
\]
\end{theorem}
However for general symplectic manifolds, the isomorphisms $H^{p}( A^{\ast}_{coE}(M))\cong \tilde H^{p}(A^{\ast}(M))$ does not hold.
In fact  counter examples are given in \cite{FIL}.
So far we have hardly found examples of non-K\"ahler manifolds such that isomorphisms $H^{p}( A^{\ast}_{coE}(M))\cong \tilde H^{p}(A^{\ast}(M))$ hold.
The purpose of this paper is to compute the coeffective cohomology of some class of symplectic manifolds by using of finite dimensional cochain complex, and give  non-K\"ahler examples such that isomorphisms $H^{p}( A^{\ast}_{coE}(M))\cong \tilde H^{p}(A^{\ast}(M))$ hold.

\section{Preliminary: Coeffective cohomology of sub-complex}
Let $(M, \omega)$ be a compact $2n$-dimensional symplectic manifold.
\begin{proposition}\label{suj}{\rm (\cite{FIL},\cite{Lie})}
Then the map $\omega\wedge: A^{p}(M)\to A^{p+2}(M)$ is injective for $p\le n-1$ and surjective for $p\ge n-1$.
\end{proposition}
By this proposition we have $H^{p}(A^{\ast}_{coE}(M))=\{0\}$ for $p\le n-1$ and so it is sufficient to consider $H^{p}(A^{\ast}_{coE}(M))$ for  $p\ge n$.
Since $\omega$ is closed, we have the short exact sequence of cochain complexes
\[\xymatrix{
	0\ar[r]&A_{coE}^{\ast}(M)\ar[r]& A^{\ast}(M)\ar[r]^{\omega \wedge} &\omega \wedge A^{\ast}(M) \ar[r]&0,
 }
\]
where we consider $\omega \wedge A^{\ast}(M) $ the cochain complex which is graded as $(\omega \wedge A^{\ast}(M))^{p}=\omega \wedge A^{p}(M)$.
By this sequence 
we have the long exact sequence of cohomology
\[\xymatrix{
	\  \ar[r]& H^{p-1}(A^{\ast}(M)) \ar[r]^{(\omega \wedge)^{\ast}}& H^{p+1}(\omega \wedge A^{\ast}(M)) \ar[r]&H^{p}(A_{coE}^{\ast}(M))\ar[r]& H^{p}(A^{\ast}(M))\ar[r]^(0,8){(\omega \wedge)^{\ast}}&\ 
 }.
\]
By Proposition \ref{suj},  we have  $\omega\wedge A^{p-1}(M)=A^{p+1}(M)$ for $p\ge n$ and so  the exact sequence is given by
\[\xymatrix{
	\ & \ar[r]& H^{p-1}(A^{\ast}(M))  \ar[r]^{(\omega \wedge)^{\ast}}& H^{p+1}(A^{\ast}(M)) \ar[r]&H^{p}(A_{coE}^{\ast}(M))\ar[r]&  H^{p}(A^{\ast}(M))\ar[r]^(0,8){(\omega \wedge)^{\ast}}&\
 }.
\]

\begin{proposition}\label{juu}
Let $A^{\ast}\subset A^{\ast}(M)$ be a sub-complex such that the inclusion $\Phi :A^{\ast}\to A^{\ast}(M)$ induces a cohomology isomorphism.
Assume $\omega\in A^{\ast}$ and the map $\omega\wedge :A^{p}\to A^{p+2}$ is surjective for $p\ge n-1$.
Denote $A^{\ast}_{coE}=\ker (\omega \wedge)_{\vert A^{\ast}}$.
Then  the inclusion $\Phi :A_{coE}^{\ast}\to A_{coE}^{\ast}(M)$ induces an isomorphism
\[H^{p}(A^{\ast}_{coE})\cong H^{p}(A^{\ast}_{coE}(M))
\]
for $p\ge n$.
\end{proposition}
\begin{proof}
As above, we have the exact sequence of cochain complex
\[\xymatrix{
	0\ar[r]&A^{\ast}_{coE}\ar[r]&A^{\ast}\ar[r]^{\omega\wedge }&\omega \wedge A^{\ast}\ar[r]&0.
 }
\]
By the assumption, for $p\ge n$ we have the long exact sequence of cohomology
\[\xymatrix{
	\  \ar[r]& H^{p-1}(A^{\ast}) \ar[r]^{(\omega \wedge)^{\ast}}& H^{p+1}(A^{\ast}) \ar[r]&H^{p}(A^{\ast}_{coE})\ar[r]&  H^{p}(A^{\ast})\ar[r]^(0,8){(\omega \wedge)^{\ast}}&\ 
 }.
\]
By the inclusion $\Phi:(\bigwedge_{coE} A ^{\ast})^{T}\to A^{\ast}_{coE}(M)$,
we have the commutative diagram
\[\xymatrix{
 H^{p-1}(A^{\ast}(M)) \ar[r]^{(\omega \wedge)^{\ast}}& H^{p+1}(A^{\ast}(M)) \ar[r]&H^{p}(A_{coE}^{\ast}(M))\ar[r]&H^{p}(A^{\ast}(M))\ar[r]^{(\omega \wedge)^{\ast}} & H^{p+2}(A^{\ast}(M)) \\
H^{p-1}(A^{\ast})\ar[u]^{\Phi^{\ast}}  \ar[r]^{(\omega \wedge)^{\ast}}& H^{p+1}(A^{\ast}) \ar[r]\ar[u]^{\Phi^{\ast}}&H^{p}(A^{\ast}_{coE})\ar[r]\ar[u]^{\Phi^{\ast}}&H^{p}(A^{\ast})   \ar[u]^{\Phi^{\ast}}\ar[r]^{(\omega \wedge)^{\ast}}&H^{p+2}(A^{\ast}) \ar[u]^{\Phi^{\ast}}.
 }
 \]
By the assumption   $\Phi^{\ast}:H^{\ast}(A^{\ast})\to H^{\ast}(A^{\ast}(M))$ is an isomorphism and so
 by this diagram $\Phi^{\ast}:H^{p}(A^{\ast}_{coE})\to H^{p}(A^{\ast}_{coE}(M))$ is an isomorphism.

\end{proof}

\section{Background: Fernandez-Ibanez-de Leon's theorem}
Let $G$ be a simply connected Lie group with a lattice (i.e. a cocompact discrete subgroup of $G$) $\Gamma$.
We call $G/\Gamma$ a nilmanifold (resp. solvmanifold) if $G$ is nilpotent (resp. solvable).
Let $\g$ be the Lie algebra of $G$ and $\bigwedge \g^{\ast}$ be the cochain complex of $\g$ with the differential which is induced by the dual of the Lie bracket.
As we regard $\bigwedge \g^{\ast}$ as the left-invariant forms on $G/\Gamma$, we consider the inclusion $\bigwedge \g^{\ast}\subset A^{\ast}(G/\Gamma)$.
Let $\omega \in \bigwedge^{2}\g^{\ast}$ be a left-invariant symplectic form.
Then the map $\omega\wedge : \bigwedge ^{p}\g^{\ast}\to  \bigwedge ^{p+2}\g^{\ast}$ is surjective for $p\ge n-1$(see \cite{FIL}).
In \cite{Nom} Nomizu showed that if $G$ is nilpotent then  the inclusion  $\bigwedge \g^{\ast}\subset A^{\ast}(G/\Gamma)$ induces an isomorphism of cohomology.
Hence by Proposition \ref{juu}, we have the following theorem which was noted in \cite{FIL} and \cite{FIL2}.
\begin{theorem}\label{Fer}
Let $G$ be a simply connected   nilpotent Lie group with a lattice $\Gamma$ and a left-invariant symplectic form $\omega$.
Then the inclusion $\bigwedge \g^{\ast}\subset A^{\ast}(G/\Gamma)$ induces an isomorphism 
\[
H^{p}(\bigwedge_{coE} \g^{\ast})\cong H^{p}(A^{\ast}_{coE}(G/\Gamma))
\]
for $p\ge n$ where $\bigwedge_{coE} \g^{\ast}=\{\alpha\in\bigwedge {\frak g}^{\ast} \vert  \omega\wedge \alpha=0\}$.
\end{theorem}
 In \cite{Hatt} Hattori showed that the isomorphism $H^{\ast}(\bigwedge \g^{\ast})\cong H^{\ast}(A^{\ast}(G/\Gamma))$ also holds if $G$ is completely solvable (i.e. $G$ is solvable and for any $g\in G$ the all eigenvalues of the adjoint operator ${\rm Ad}_{g}$ are real).
 Thus  we can extend this theorem for completely solvmanifolds.
 However for a general solvmanifold $G/\Gamma$, the isomorphism $H^{\ast}(\bigwedge_{coE} \g^{\ast})\cong H^{\ast}(A^{\ast}_{coE}(G/\Gamma))$ does not hold and we can't compute the coeffective cohomology by using of $\bigwedge \g^{\ast}$.

In \cite{B} Baues constructed compact aspherical manifolds $M_{\Gamma}$ such that  the class of these aspherical manifolds contains the class of solvmanifolds and showed that the de Rham cohomology of these aspherical manifolds can be computed by certain finite dimensional cochain complexes. 
In next section by using of Baues's results, we will show a generalization of  Theorem \ref{Fer}.

\section{Main results: Coeffective cohomology of aspherical  manifolds with torsion-free virtually polycyclic fundamental groups} 

\subsection{Notation and conventions}
Let $k$ be a subfield of $\C$.
A group $\bf G$ is called a $k$-algebraic group if $\bf G$ is a Zariski-closed subgroup of $GL_{n}(\C)$ which is defined by polynomials with coefficients in $k$.
Let  ${\bf G}(k)$ denote the set  of  $k$-points of $\bf G$ 
and ${\bf U}({\bf G})$ the maximal Zariski-closed unipotent normal $k$-subgroup of $\bf G$ called the unipotent radical of $\bf G$.
If $\bf G$ consists of  semi-simple elements, we call $\bf G$ a d-group.
Let $U_{n}(k)$ denote the $n\times n$ $k$-valued upper triangular unipotent matrix group.

\subsection{Baues's results}
A group $\Gamma$ is called polycyclic if it admits a sequence 
\[\Gamma=\Gamma_{0}\supset \Gamma_{1}\supset \cdot \cdot \cdot \supset \Gamma_{k}=\{ e \}\]
of subgroups such that each $\Gamma_{i}$ is normal in $\Gamma_{i-1}$ and $\Gamma_{i-1}/\Gamma_{i}$ is cyclic.
We denote ${\rm rank}\,\Gamma=\sum_{i=1}^{i=k}{\rm rank}\, \Gamma_{i-1}/\Gamma_{i}$.
We define an infra-solvmanifold as a manifold of the form $G/ \Delta$ where $G$ is a simply connected solvable Lie group, and $\Delta $ is a torsion free subgroup of ${\rm Aut}(G)\ltimes G$ such that for the projection $p: {\rm Aut}(G)\ltimes G\to  {\rm Aut}(G)$ $p(\Delta)$ is contained in a compact subgroup of ${\rm Aut}(G)$.
By a result of Mostow in \cite{Mos}, the fundamental group of an infra-solvmanifold is virtually polycyclic(i.e it contains a finite index polycyclic subgroup).
In particular, a lattice $\Gamma$ of a simply connected solvable Lie group $G$ is a polycyclic group with ${\rm rank}\,  \Gamma=\dim G$(see \cite{R}). 

Let $k$ be a subfield of $\mathbb{C}$.
Let $\Gamma$ be a torsion-free virtually polycyclic group.
For a finite index polycyclic subgroup $\Delta\subset \Gamma$, we denote ${\rm rank}\, \Gamma ={\rm  rank}\, \Delta$.

\begin{definition}\label{a-d-2}
We call a $k$-algebraic group  ${\bf H}_{\Gamma}$ a $k$-algebraic hull of $\Gamma$ if there exists an injective group homomorphism $\psi:\Gamma\to {\bf H}_{\Gamma}(k)$
and ${\bf H}_{\Gamma}$ satisfies the following conditions:
\\
(1)  \ $\psi (\Gamma)$ is Zariski-dense in $\bf{H}_{\Gamma}$.\\
(2) \   $Z_{{\bf H}_{\Gamma}}({\bf U}({\bf H}_{\Gamma}))\subset {\bf U}({\bf H}_{\Gamma})$, where  $Z_{{\bf H}_{\Gamma}}({\bf U}({\bf H}_{\Gamma}))$ is the centralizer of ${\bf U}({\bf H}_{\Gamma})$.\\
(3) \ $\dim {\bf U}({\bf H}_{\Gamma})$=${\rm rank}\,\Gamma$.   
\end{definition}
\begin{theorem}\label{a-t-1}{\rm (\cite[Theorem A.1]{B})}
There exists a $k$-algebraic hull of $\Gamma$ and a k-algebraic hull of $\Gamma$ is unique up to $k$-algebraic group isomorphism.
\end{theorem}

 Let $\Gamma$ be a torsion-free virtually polycyclic group and
$\bf{H_{\Gamma}}$  the $\Ra$-algebraic hull of $\Gamma$.
Denote $H_{\Gamma}=\bf{H_{\Gamma}}(\R)$. 
Let $U_{\Gamma}$ be the unipotent radical of $H_{\Gamma}$ and $T$  a maximal  d-subgroup.
Then $H_{\Gamma}$ decomposes as a semi-direct product $H_{\Gamma}=T\ltimes U_{\Gamma}$ see cite[Proposition 2.1]{B}.
Let $\frak{u}$ be the Lie algebra of $U_{\Gamma}$. Since the exponential map ${\exp}:{\frak u} \longrightarrow U_{\Gamma}$ is a diffeomorphism, $U_{\Gamma}$ is diffeomorphic to $\R^n$ such that $n={\rm rank}\,\Gamma$.
For the semi-direct product $H_{\Gamma}=T\ltimes U_{\Gamma}$, we denote $\phi:T\to {\rm Aut}(U_{\Gamma})$ the action of $T$ on $U_{\Gamma}$.
Then we have the homomorphism $\alpha :H_{\Gamma}\longrightarrow {\rm Aut}(U_{\Gamma})\ltimes U_{\Gamma}$ such that $\alpha(t,u)=(\phi(t), u)$ for $(t,u)\in T\ltimes U_{\Gamma}$.
By the property (2) in Definition \ref{a-d-2}, $\phi$ is injective and hence $\alpha$ is injective.

In \cite{B} Baues constructed a compact  aspherical manifold $M_{\Gamma}=\alpha(\Gamma)\backslash U_{\Gamma}$ with $\pi_{1}(M_{\Gamma})=\Gamma$.
We call $M_{\Gamma}$ a standard $\Gamma$-manifold.

\begin{theorem}{\rm (\cite[Theorem 1.2, 1.4]{B})}\label{a-t-3}
A standard $\Gamma$-manifold is unique up to diffeomorphism.
A compact infra-solvmanifold with the fundamental group $\Gamma$ is diffeomorphic to  the standard $\Gamma$-manifold $M_{\Gamma}$.
In particular, a solvmanifold  $G/\Gamma$ is diffeomorphic to the standard $\Gamma$-manifold $M_{\Gamma}$.
\end{theorem}
Let $A^{\ast}(M_{\Gamma})$ be the de Rham complex of $M_{\Gamma}$.
Then $A^{\ast}(M_{\Gamma}) $ is  the set of   the $\Gamma$-invariant differential forms ${A^{\ast}(U_{\Gamma})}^{\Gamma}$ on $U_{\Gamma}$. 
Let $(\bigwedge \frak{u} ^{\ast})^{T}$ be the left-invariant forms on $U_{\Gamma}$ which are fixed by $T$.
Since $\Gamma\subset H_{\Gamma}=U_{\Gamma}\cdot T$, we have the inclusion
\[(\bigwedge {\frak u} ^{\ast})^{T} ={A^{\ast}(U_{\Gamma})}^{H_{\Gamma}} \subset {A^{\ast}(U_{\Gamma})}^{\Gamma}= A^{\ast}(M_{\Gamma}).\]
\begin{theorem}{\rm (\cite[Theorem 1.8]{B})}\label{a-t-4}
This inclusion induces a cohomology isomorphism.
\end{theorem}
\subsection{Main results}
Let $\omega\in (\bigwedge {\frak u} ^{\ast})^{T}$ be a symplectic form.
Denote 
\[\bigwedge_{coE} {\frak u} ^{\ast}=\{\alpha\in\bigwedge {\frak u}^{\ast} \vert \omega\wedge \alpha=0\}\]
and
\[\tilde H^{\ast}((\bigwedge {\frak u} ^{\ast})^{T})=\{[\alpha]\in H^{\ast}((\bigwedge {\frak u} ^{\ast})^{T})\vert [\omega]\wedge [\alpha]=0\}.
\]
By Theorem \ref{a-t-4}, we have $\tilde H^{\ast}((\bigwedge {\frak u} ^{\ast})^{T}) \cong \tilde H^{\ast}(A^{\ast}(M_{\Gamma}))$.

\begin{lemma}\label{keyl}
For $p\ge n-1$, the linear map $\omega \wedge:(\bigwedge^{p}{\frak u} ^{\ast})^{T}\to (\bigwedge ^{p+2}{\frak u} ^{\ast})^{T}$ is surjective.
\end{lemma}
\begin{proof}
First we notice that the map $\omega \wedge:\bigwedge^{p} {\frak u} ^{\ast}\to \bigwedge ^{p+2}{\frak u} ^{\ast}$ is surjective (see \cite[Lemma 2.1]{FIL}).
Since $T$ is d-group, for $t\in T$ the $t$-action on $\bigwedge {\frak u} ^{\ast}$ is diagonalizable (see \cite{B}).
Hence we have a decomposition
\[\bigwedge^{p}{\frak u} ^{\ast}=A^{p}\oplus B^{p}\]
such that $A^{p}$ is the subspace of $t$-invariant elements and $B^{p}$ is its complement.
Since the $t$-action is diagonalizable, we have a basis $\{ x_{1},\dots x_{2n}\}$ of $ {\frak u}^{\ast}\otimes \C$ such that the $t$-action is represented by a diagonal matrix.
Then we have
\[A^{p}\otimes \C=\langle x_{i_{1}}\wedge\dots\wedge x_{i_{p}}\vert 1\le i_{1}<\dots<i_{p}\le 2n,\, t\cdot( x_{i_{1}}\wedge\dots \wedge x_{i_{p}})=x_{i_{1}}\wedge\dots \wedge x_{i_{p}}\rangle,
\]
and
\begin{multline*}
B^{p}\otimes \C=\\
\langle x_{i_{1}}\wedge\dots \wedge x_{i_{p}}\vert 1\le i_{1}<\dots<i_{p}\le 2n,\, t\cdot ( x_{i_{1}}\wedge\dots \wedge x_{i_{p}})=\alpha_{i_{1}\dots i_{p}}(t)x_{i_{1}}\wedge\dots \wedge x_{i_{p}} , \alpha_{i_{1}\dots i_{p}}(t)\not=1 \rangle.
\end{multline*}
By $\omega\in (\bigwedge {\frak u} ^{\ast})^{T}$, we have $\omega=\sum  a_{kl}x_{k}\wedge x_{l}$ such that if $a_{kl}\not =0$, then $x_{k}\wedge x_{l}\in A^{p}\otimes \C$.
Then for $x_{i_{1}}\wedge \dots \wedge x_{i_{p}}\in B^{p}\otimes \C$   we have
 \[\omega \wedge x_{i_{1}}\wedge \dots \wedge x_{i_{p}}=\sum a_{kl}x_{k}\wedge x_{l}\wedge x_{i_{1}}\wedge \dots \wedge x_{i_{p}}.\]
If $a_{kl}\not =0$, we have
\[t\cdot (x_{k}\wedge x_{l}\wedge x_{i_{1}}\wedge \dots \wedge x_{i_{p}})=\alpha_{i_{1}\dots i_{p}}(t)x_{k}\wedge x_{l}\wedge x_{i_{1}}\wedge \dots \wedge x_{i_{p}}.
\]
Thus $\omega \wedge x_{i_{1}}\wedge \dots \wedge x_{i_{p}}\in B^{p+2}\otimes \C$.
By this we have $(\omega\wedge B^{p})\subset B^{p+2}$.
Since $T$ acts semi-simply on $\bigwedge ^{p}{\frak u} ^{\ast}$,  we consider the decomposition 
\[\bigwedge ^{p}{\frak u} ^{\ast}=(\bigwedge ^{p}{\frak u} ^{\ast})^{T}\oplus C^{p}\]
such that $C^{p}$ is a complement of $(\bigwedge ^{p}{\frak u} ^{\ast})^{T}$ for $T$-action.
By the above argument we have $(\omega \wedge C^{p})\subset C^{p+2}$.
Clearly we have $(\omega \wedge (\bigwedge ^{p}{\frak u} ^{\ast})^{T})\subset (\bigwedge ^{p+2}{\frak u} ^{\ast})^{T}$.
Since for $p\ge n-1$ the map $\omega \wedge:\bigwedge^{p} {\frak u} ^{\ast}\to \bigwedge ^{p+2}{\frak u} ^{\ast}$ is surjective, we have 
\[(\bigwedge ^{p+2}{\frak u} ^{\ast})^{T}\oplus C^{p}=\omega \wedge  \bigwedge ^{p}{\frak u} ^{\ast}=(\omega \wedge  (\bigwedge ^{p}{\frak u} ^{\ast})^{T} )\oplus (\omega \wedge C^{p}).\]
Thus we have $\omega \wedge  (\bigwedge ^{p}{\frak u} ^{\ast})^{T} =  (\bigwedge ^{p+2}{\frak u} ^{\ast})^{T} $.
Hence the lemma follows.
\end{proof} 

By this lemma and Proposition \ref{juu}, we have:
\begin{theorem}\label{MMMMT}
Let $\Gamma$ be a torsion-free virtually polycyclic group and $M_{\Gamma}$ the standard $\Gamma$-manifold with
 a symplectic form $\omega$ such that $\omega\in (\bigwedge {\frak u} ^{\ast})^{T}$.
Then for $p\ge n$, the inclusion $\Phi:(\bigwedge_{coE} {\frak u} ^{\ast})^{T}\to A^{\ast}_{coE}(M_{\Gamma})$
  induces  an isomorphism $\Phi^{\ast}:H^{\ast}((\bigwedge_{coE} {\frak u} ^{\ast})^{T})\cong H^{\ast}(A^{p}_{coE}(M_{\Gamma}))$.
\end{theorem}

\begin{remark}\label{CS}
In \cite{K2}, the author showed that  if there exists $[\omega]\in H^{2}(M_{\Gamma},\R)$ such that $[\omega]^{\frac{1}{2}\dim M_{\Gamma}}\not=0$, then an invariant form $\omega\in (\bigwedge {\frak u} ^{\ast})^{T}$ which represent the cohomology class $[\omega]$ is a symplectic form on $M_{\Gamma}$.
Hence if $M_{\Gamma}$ is cohomologically symplectic(i.e. there exists $[\omega]\in H^{2}(M_{\Gamma},\R)$ such that $[\omega]^{\frac{1}{2}\dim M_{\Gamma}}\not=0$),  then $M_{\Gamma}$ admits a symplectic form $\omega$ such that $\omega\in (\bigwedge {\frak u} ^{\ast})^{T}$.
\end{remark}
\begin{corollary}\label{abb}
Under the same assumption of Theorem \ref{MMMMT},
if $U_{\Gamma}$ is abelian, then for $p\ge n$ we have an isomorphism
\[ H^{p}( A^{\ast}_{coE}(M_{\Gamma}))\cong \tilde H^{p}(A^{\ast}(M_{\Gamma})).\]
\end{corollary}
\begin{proof}
If  $U_{\Gamma}$ is abelian, then the differential of $\bigwedge {\frak u} ^{\ast}$ is $0$.
Hence we have 
\[H^{\ast}(A^{\ast}(M_{\Gamma}))\cong H^{\ast}((\bigwedge {\frak u} ^{\ast})^{T})= (\bigwedge {\frak u} ^{\ast})^{T}
\]
and 
\[H^{\ast}((\bigwedge_{coE} {\frak u} ^{\ast})^{T})= (\bigwedge_{coE} {\frak u} ^{\ast})^{T}.\]
This gives 
\[\tilde H^{\ast}(A^{\ast}(M_{\Gamma})) \cong\tilde H^{\ast}((\bigwedge {\frak u} ^{\ast})^{T})=\{\alpha\in (\bigwedge {\frak u} ^{\ast})^{T}\vert \alpha\wedge \omega=0\}=(\bigwedge_{coE} {\frak u} ^{\ast})^{T}=H^{\ast}((\bigwedge_{coE} {\frak u} ^{\ast})^{T}).
\]
Hence by the above theorem the corollary follows.
\end{proof}
In \cite{K} the author showed the following theorem.

\begin{theorem}\label{thoo}{\rm (\cite{K})}
Let $\Gamma$ be a torsion-free virtually polycyclic group.
Then the following two conditions are equivalent:\\
$(1)$ ${\bf U}_{\Gamma}$ is abelian.\\
$(2)$ $\Gamma$  is a finite extension group of a lattice of a Lie group  $G=\R^{n}\ltimes_{\phi} \R^{m}$ such that the action $\phi:\R^{n}\to {\rm  Aut}(\R^{m})$ is semi-simple.
\end{theorem}
Hence we have:
\begin{corollary}\label{CCC}
Under the same assumption of Theorem \ref{MMMMT},
if  $\Gamma$ satisfies the condition $(2)$ in Theorem \ref{thoo},
then for $p\ge n$ we have an isomorphism 
\[
 H^{p}( A^{\ast}_{coE}(M_{\Gamma}))\cong \tilde H^{p}(A^{\ast}(M_{\Gamma})).\]
\end{corollary}
\begin{remark}
In fact by Arapura and Nori's theorem(\cite{AN}) a virtually polycyclic group $\Gamma$ must be virtually abelian if the standard $\Gamma$-manifold is K\"ahler. Therefore $G/\Gamma$ is finitely covered by a torus and the assumptions of \ref{thoo} are satisfied.
By Arapura and Nori's theorem, if a solvmanifold $G/\Gamma$ admits  a K\"ahler structure, then $G$ is (I)-type (i.e. for any $g\in G$  all eigenvalues of the adjoint operator ${\rm Ad}_{g}$ have absolute value 1).
Thus in the above corollary if $G$ is not (I)-type, then $M_{\Gamma}$ does not admit a K\"ahler structure.
The author gave such non-K\"ahler examples in \cite{K}. 
\end{remark}

\section{examples}

\begin{example}
First we give examples of solvmanifolds such that  $H^{p}( A^{\ast}_{coE}(M_{\Gamma}))\cong \tilde H^{p}(A^{\ast}(M_{\Gamma}))$ by using of Corollary \ref{CCC}.
We notice that if a solvmanifold $G/\Gamma$ has a symplectic form $\omega$ then we have a closed two form $\omega_{0}\in (\bigwedge {\frak u} ^{\ast})^{T}$ which is homologous to $\omega$ and $\omega_{0}$ is also a symplectic form as we note in Remark \ref{CS}.
Let $G=\C\ltimes_{\phi} \C^{2}$ with $\phi(x)= \left(
\begin{array}{ccc}
e^{x}&0\\
0&e^{-x}
\end{array}
\right)$.
Then it is known that $G$ has a left-invariant symplectic  form and a lattice $\Gamma$ (see \cite{Na}).
Thus we have a symplectic form $\omega\in  (\bigwedge {\frak u} ^{\ast})^{T}$ and by Corollary \ref{CCC} we have an isomorphism $H^{p}( A^{\ast}_{coE}(G/\Gamma))\cong \tilde H^{p}(A^{\ast}(G/\Gamma))$.
\begin{remark}
$G$ is not completely solvable.
In fact the de Rham cohomology of $G/\Gamma$ varies  according to a choice of a lattice $\Gamma$.
Thus it is not easy to compute the coeffective cohomology of $G/\Gamma$ by using of $\bigwedge \g^{\ast}$.
\end{remark}
\begin{remark}
$G$ is not (I)-type and hence $G/\Gamma$ does not admit a K\"ahler structure. 
\end{remark}
\end{example}

\begin{example}
We give an example of a symplectic manifold $M_{\Gamma}$ such that the isomorphism  $H^{p}( A^{\ast}_{coE}(M_{\Gamma}))\cong \tilde H^{p}(A^{\ast}(M_{\Gamma}))$ holds but $U_{\Gamma}$ is not abelian.
Let $\Gamma=\Z \ltimes _{\phi}\Z^{2}$ such that 
 for $t\in\Z$
\[\phi(t)
=\left(
\begin{array}{cc}
(-1)^{t}& (-1)^{t} t  \\
0&     (-1)^{t}  
\end{array}
\right).
\]
Then we have $H_{\Gamma}= \{\pm 1\} \ltimes  U_{3}(\R)$ such that

\[ (-1)\cdot \left(
\begin{array}{ccc}
1&  x&  z  \\
0&     1&y      \\
0& 0&1 
\end{array}
\right)=
\left(
\begin{array}{ccc}
1&  x&  (-1)z  \\
0&     1&(-1)y      \\
0& 0&1 
\end{array}
\right)
\]
(see \cite[Section 7]{K}).
The dual space of  the Lie algebra $\frak u$ of $U_{\Gamma}$ is given by 
${\frak u}^{\ast}=\langle x_{1}, x_{2}, x_{3}\rangle$ such that the 
differential is given by 
\[dx_{1}=dx_{2}=0,\  dx_{3}=-x_{1}\wedge x_{2}.\]
The action of $\{\pm 1 \}$ on $U_{\Gamma}$ is given by
\[(-1)\cdot x_{1}=x_{1},
\] 
\[(-1)\cdot x_{2}=-x_{2},\, (-1)\cdot x_{3}=-x_{3}.\]
Then we have $(\bigwedge {\frak u}^{\ast})^{\{\pm 1\}}=\bigwedge \langle x_{1}, x_{2}\wedge x_{3}\rangle$.
By this the differential on $(\bigwedge {\frak u}^{\ast})^{\{\pm 1\}}$ 
is $0$.
We consider the product $M_{\Gamma}\times M_{\Gamma}$ for this $\Gamma$.
Then by the cochain complex $(\bigwedge {\frak u}^{\ast})^{\{\pm 1\}}\otimes (\bigwedge {\frak u}^{\ast})^{\{\pm 1\}}=\bigwedge \langle x_{1}, 
x_{2}\wedge x_{3}\rangle\otimes \bigwedge \langle y_{1}, y_{2}\wedge y_{3}\rangle$ we can compute the de Rham cohomology and coeffective 
cohomology of $M_{\Gamma}\times M_{\Gamma}$ where we denote $y_{1}, y_{2}, y_{3}$ the copy of $x_{1}, x_{2}, x_{3}$.
We have a symplectic form 
\[\omega=x_{1}\wedge y_{1}+x_{2}\wedge x_{3}+y_{2}\wedge y_{3}\]
on $M_{\Gamma}\times M_{\Gamma}$.
Then we have:
\begin{proposition}
 For $p\ge n$   we have an isomorphism 
\[H^{p}( A^{\ast}_{coE}(M_{\Gamma}\times M_{\Gamma} ))\cong \tilde H^{p}(A^{\ast}(M_{\Gamma}\times M_{\Gamma} )).
\]
\end{proposition}
\begin{proof}
Since the differential on $(\bigwedge {\frak u}^{\ast})^{\{\pm 1\}}\otimes (\bigwedge {\frak u}^{\ast})^{\{\pm 1\}}$ is not $0$ as above, the 
proposition follows as the proof of Corollary \ref{abb}.
\end{proof}
\begin{remark}
$M_{\Gamma}$ is finitely covered by a quotient of $U_{3}(\R)$ by a 
lattice.
Thus $M_{\Gamma}\times M_{\Gamma}$ is finitely covered by the product of 
such nilmanifolds.
The de Rham cohomology and coeffective cohomology of this covering space 
are computed by $\bigwedge {\frak u}^{\ast}\otimes \bigwedge {\frak u}^{\ast}$.
This space does not satisfy the isomorphism as this proposition.
Indeed $x_{1}\wedge x_{2}\wedge y_{2}\wedge y_{3}$ is coeffective and 
its coeffective cohomology class is not $0$.
But we have $d(x_{3}\wedge y_{2}\wedge y_{3})=x_{1}\wedge x_{2}\wedge y_{2}\wedge y_{3}$ and hence its de Rham cohomology class  is $0$.
Thus we have 
\[H^{4}( A^{\ast}_{coE}((U_{3}(\R)/\Gamma^{\prime})\times (U_{3}(\R)/\Gamma^{\prime})))\not\cong \tilde H^{4}(A^{\ast}((U_{3}(\R)/\Gamma^{\prime})\times (U_{3}(\R)/\Gamma^{\prime}))).\]

\end{remark}

\end{example}

{\bf  Acknowledgements.} 

The author would like to express his gratitude to   Toshitake Kohno for helpful suggestions and stimulating discussions.
This research is supported by JSPS Research Fellowships for Young Scientists.

\end{document}